\documentclass[12pt]{article}
\usepackage{amsfonts, amsmath, latexsym, epic, eepic, pifont, url}
\title{Infinite serie of extreme Delaunay polytope}
\author{Mathieu DUTOUR\\
        \normalsize  LIGA, ENS/CNRS, Paris and\\
        \normalsize  Hebrew University, Jerusalem\footnote{Research financed by 
EC's IHRP Programme, within the Research Training Network ``Algebraic Combinator
ics in Europe,'' grant HPRN-CT-2001-00272.}\\
}

\begin{document}
\newcommand{\R}{\ensuremath{\mathbb{R}}}
\newcommand{\N}{\ensuremath{\mathbb{N}}}
\newcommand{\Q}{\ensuremath{\mathbb{Q}}}
\newcommand{\C}{\ensuremath{\mathbb{C}}}
\newcommand{\Z}{\ensuremath{\mathbb{Z}}}
\newcommand{\T}{\ensuremath{\mathbb{T}}}
\newtheorem{prop}{Proposition}
\newtheorem{theorem}{Theorem}
\newtheorem{corollary}{Corollary}
\newtheorem{lem}{Lemma}
\newtheorem{conjecture}{Conjecture}
\newtheorem{claim}{Claim}
\newtheorem{remark}{Remark}
\newtheorem{definition}{Definition}
\newcommand{\qed}{\hfill $\Box$ }
\newcommand{\proof}{\noindent{\bf Proof.}\ \ }

\maketitle

\begin{abstract}
\noindent A Delaunay polytope $P$ is said to be {\em extreme} if the only (up to isometries) affine bijective transformations $f$ of $\R^n$, for which $f(P)$ is again a Delaunay polytope, are the homotheties. 
This notion was introduced in \cite{DGL92}; also some examples in dimension $1$, $6$, $7$, $15$, $16$, $22$, $23$ were constructed and it was proved that in dimension less than $6$ there are no extreme Delaunay polytopes, except the segment.
In this note, for every $n\geq 6$ we build an extreme Delaunay polytope $ED_n$ of dimension $n$. 
\end{abstract}

Let $L\subset \R^n$ be a $n$-dimensional lattice and let $S=S(c,r)$ be a  sphere in $\R^n$ with center $c$ and radius $r$. Then, $S$ is said to be an {\em empty sphere} in $L$ if the following two conditions hold:
\begin{center}
$\Vert v-c\Vert\geq r$ for all $v\in L$ and the set $S\cap L$ has affine rank $n+1$.
\end{center}
The polytope $P$, which is defined as the convex hull of the set $P=S\cap L$,
is called a {\em Delaunay polytope}.

Denote by $D_n$ the root lattice defined as
\begin{equation*}
D_n=\{(x_1, \dots, x_n)\in\Z^n\mbox{~with~}\sum_{i=1}^n x_i\mbox{~even}\}\;.
\end{equation*}
The Delaunay polytopes of this lattice are:
\begin{itemize}
\item {\em Half-cube} $\frac{1}{2}H_n=\{x\in \{0,1\}^n\mbox{~with~}\sum_{i=1}^n x_i\mbox{~even}\}$; it has $2^{n-1}$ vertices.
\item {\em Cross-polytope} $C_n=\{e_1\pm e_i\mbox{~with~}1\leq i\leq n\}$ with $e_i=(0,\dots, 1,\dots,0)$; it has $2n$ vertices.
\end{itemize}

Define the polytope $P_n$ to be the convex hull of the three following layers of points:
\begin{itemize}
\item a vertex $V=(\frac{1}{2},\dots, \frac{1}{2},1)$,
\item $2^{n-2}$ vertices $(x_1,\dots, x_{n-1},0)$ with $x_i\in\{0,1\}$ and $\sum_{i=1}^{n-1} x_i$ even,
\item $2(n-1)$ vertices $V_{j,\pm}=(\frac{1}{2},\dots, \frac{1}{2},-1)\pm e_j$ with $1\leq j\leq n-1$.
\end{itemize}

It is easy to see, that, denoting by $L_n$ the lattice generated by all $(x,0)$ with $x\in D_{n-1}$ and $(\frac{1}{2},\dots, \frac{1}{2},1)$, the third layer belongs to $L_n$ if and only if $n$ is even.

\begin{theorem}\label{TheoremInfiniteSeq}
Let $n$ be even and greater than $6$.

(i) $P_n$ is an extreme Delaunay polytope; the corresponding quadratic form is $q(x)=x_1^2+\dots+x_{n-1}^2+\frac{n-3}{4}x_n^2$.

(ii) The center of the circumscribing sphere is $(0,\dots,0,\frac{-1}{n-3})$, it radius is $\frac{n-2}{\sqrt{n-3}}$. 

(iii) $P_6$=Schl\"afli polytope; its group of isometries has size $51840$ and is transitive on vertices.
If $n>6$, then its group of isometries is $Sym(n-1)\times 2^{n-2}$ and there are three orbits of vertices.
\end{theorem}
\proof In order to prove extremality, we will use geometry of numbers, i.e. in order to show that the only transformations preserving the property of being Delaunay are modulo isometries and translations the homotheties, we will determine the set of quadratic forms $q$ such that $P_n$ is a Delaunay polytope for $q$.

The restriction of the quadratic form $q$ to the plane $x_n=0$ left us with the $2^{n-2}$ vertices $(x_1, \dots, x_{n-1}, 0)$ with $x_i\in \{0,1\}$ and $\sum_{i} x_i$ even.
The rank of Half-cube is $n$ (see \cite{DL} and \cite{Du}), we can write $q$ restricted to the plane $x_n=0$ as $\sum_{i=1}^{n-1} \alpha_i x_i^2$ with $\alpha_i>0$. We write the form $q$ in the following way:
\begin{equation*}
q(x)=\sum_{i=1}^{n-1}\alpha_i (x_i+\beta_i x_n)^2+\gamma x_n^2\;.
\end{equation*}
Denote by $S(c,r)$ the empty sphere of center $c$ and radius $r$ circumscribing $\frac{1}{2}H_{n-1}$ and $V$. Denoting $c=(h_1, \dots, h_{n-1}, h_n)$, one obtains the equations:
\begin{equation*}
\left\lbrace\begin{array}{l}
\sum_{i=1}^{n-1}\alpha_i (h_i-x_i+\beta_i h_n)^2+\gamma h_n^2=r^2,\\
\sum_{i=1}^{n-1}\alpha_i (h_i-\frac{1}{2}+\beta_i (h_n-1))^2+\gamma (h_n-1)^2=r^2\;.
\end{array}\right.
\end{equation*}
The first equation is satisfied for all vectors $x_i$ with $\sum_{i=1}^{n-1}x_i$ even and $x_i\in\{0,1\}$. This implies, since $n\geq 6$ the relation $h_i+\beta_i h_n=\frac{1}{2}$. So, the equations reduce to:
\begin{equation*}
\left\lbrace\begin{array}{l}
\frac{1}{4}\sum_{i=1}^{n-1}\alpha_i+\gamma h_n^2=r^2,\\
\frac{n-2}{2}\sum_{i=1}^{n-1}\alpha_i \beta_i^2+\gamma (h_n-1)^2=r^2\;.
\end{array}\right.
\end{equation*}
Consider now the $2(n-1)$ vertices $V_{j,\pm}=(\frac{1}{2},\dots, \frac{1}{2},-1)\pm e_j$ with $1\leq j\leq n-1$.

Since $V_{j, +}$, and $V_{j, -}$ are on the sphere $S(c, r)$, one obtains:
\begin{equation*}
\left\lbrace\begin{array}{c}
\alpha_j(1-\beta_j)^2+\sum_{i\not=j} \alpha_i\beta_i^2+\gamma (h_n+1)^2=r^2,\\
\alpha_j(1+\beta_j)^2+\sum_{i\not=j} \alpha_i\beta_i^2+\gamma (h_n+1)^2=r^2\;.
\end{array}\right.
\end{equation*}
Above two equations yield $\beta_j=0$. If all $V_{j, \pm}$ are on the sphere, then $\alpha_j=\alpha$ is independent of $j$ and all equations simplify to:
\begin{equation*}
\alpha+\gamma (h_n+1)^2=r^2\mbox{~,~~}\frac{n-1}{4}\alpha+\gamma h_n^2=r^2\mbox{~~and~~}\gamma (h_n-1)^2=r^2\;,
\end{equation*}
whose solution is: 
\begin{equation*}
h_n=\frac{1}{3-n}, \;\;\gamma=\alpha\frac{n-3}{4},\;\; r=\frac{n-2}{\sqrt{n-3}}\sqrt{\alpha}\;\;. 
\end{equation*}
So, $q$ is equal to $\alpha(\sum_{i=1}^{n-1} x_i^2+\frac{n-3}{4} x_n^2)$.

Now, let us prove that this is a Delaunay polytope. If $v$ is an interior point of $P_{n}$, then it belongs to one of the three layers $x_n=-1$, $0$, $1$.
But each one of the corresponding sections (i.e. Cross-polytope, Half-cube or point) is a Delaunay polytope; so, there is no such $v$.

Let us find the symmetry group. 
Pairwise distance of adjacent vertices of Half-cube is $\sqrt{2}$. Distance between $V$ and vertices of Half-cube is $\sqrt{\frac{n-2}{2}}$. So, if $n>6$, then any isometry preserving $P_n$ must leave the vertex $V$ invariant; hence, if $n>6$, then the symmetry group is that of Half-cube, i.e. $Sym(n-1)\times 2^{n-2}$.

Standart computation shows that $P_6$ is the Schl\"afli polytope, whose symmetry group has size $51840$. \qed

Define $ED_n$ to be $P_n$ if $n$ is even and the polytope obtained from $P_{n-1}$ by using Lemma 15.3.7 is $f$ is odd. Remark that $ED_6$ (respectively $ED_7$) is Schl\"afli (respectively, Gosset) polytope.

\begin{corollary}
If $n\geq 6$, then the polytope $ED_n$ is extreme.
\end{corollary}
\proof If $n$ is even, this follows from \ref{TheoremInfiniteSeq}, while if $n$ is odd, this follows from Lemma 15.3.7. \qed

Remark, that $ED_6$ is {\em unique} extreme Delaunay polytope in dimension $6$ 
(see \cite{Hyp7}); {\em all} $6241$ 6-dimensional Delaunay polytopes were listed in \cite{Du}.


\begin{thebibliography}{99}

\bibitem[DeLa97]{DL}
M.Deza and M.Laurent, {\em Geometry of cuts and metrics}, Springer--Verlag,
 Berlin 1997.

\bibitem[Du]{Du}
M.Dutour, {\em The six-dimensional Delaunay polytopes}, submitted, \url{http://www.arxiv.org/abs/math.MG/0212353}.

\bibitem[DeDu]{Hyp7}
M.Deza and M.Dutour, {\em The hypermetric cone on seven vertices}, submitted, \url{http://www.arxiv.org/abs/math.MG/0108177}.

\bibitem[DGL92]{DGL92}
M. Deza, V.P. Grishukhin, and M. Laurent, {\em Extreme hypermetrics and L-polytopes}, in G.Hal\'asz et al. eds {\em Sets, Graphs and Numbers, Budapest (Hungary), 1991}, {\bf 60} {\em Colloquia Mathematica Societatis J\'anos Bolyai}, (1992) 157--209.




\end{thebibliography}
\end{document}